\documentclass[a4paper]{article}
\usepackage{geometry}
\usepackage{amsmath,amsfonts,amsthm}

\newcommand{\ol}[1]{\overline{#1}}

\newcommand{\II}{\mathcal{I}}

\newcommand{\eps}{\varepsilon}

\makeatletter
\newcommand{\eqcolon}{\mathrel{\mathord{=}\raise.2\p@\hbox{:}}}
\newcommand{\coloneq}{\mathrel{\raise.2\p@\hbox{:}\mathord{=}}}
\makeatother
\newcommand{\der}{\delta}

\newcommand{\RR}{\mathbb{R}}
\newcommand{\NN}{\mathbb{N}}

\newcommand{\FF}{\mathcal{F}}

\newcommand{\CC}{\mathcal{C}}
\newcommand{\BT}{\mathcal{BT}}
\newcommand{\BB}{\mathcal{B}}
\newcommand{\DD}{\mathcal{D}}
\newcommand{\ZZ}{\mathcal{Z}}
\newcommand{\LL}{\mathcal{L}}

\newcommand{\bTT}{\mathbb T}

\newcommand{\cl}{\mathcal L}

\newcommand{\crr}{\mathcal R}

\newtheorem{theorem}{Theorem}[section]

\newcommand{\TT}{\mathcal{T}}
\newcommand{\cA}{\mathcal{A}}

\usepackage{ifpdf}
\ifpdf%
  \usepackage{pdftricks}
  \begin{psinputs}
    \usepackage{pstricks}
\usepackage{pst-node}
\usepackage{pst-tree}
  \end{psinputs}
\else
  \usepackage{pstricks}
\usepackage{pst-node}
\usepackage{pst-tree}
  
\fi

\newcommand{\tnode}{\TR{\raisebox{0.5pt}{\ensuremath{\bullet}}}} 
\newcommand{\tnodel}[1]{\TR{\makebox[6pt][l]{\raisebox{0.5pt}{\ensuremath{\bullet{\text{\footnotesize
            #1}}}}}}}

\newcommand{\troot}{\bullet}

\newcommand{\tsnode}{\TR{\raisebox{0.5pt}{\ensuremath{\scriptstyle\bullet}}}} 
\newcommand{\tsroot}{\ensuremath{\scriptstyle\bullet}}
\newcommand{\aaabbb}{\pstree{\tsnode}{\pstree{\tsnode}{\tsnode}}} 
 
\newcommand{\aabb}{\pstree{\tsnode}{\tsnode}} 
\newcommand{\aababb}{\pstree{\tsnode}{\tsnode \tsnode}}

\newcommand{\smalltrees}{\psset{levelsep=-5pt,nodesep=-2pt,treesep=1pt}}
\newcommand{\largetrees}{\psset{levelsep=-10pt,nodesep=-4pt,treesep=5pt}
}
\largetrees

\newcommand{\ZCC}{\ZZ\CC}
\newcommand{\BCC}{\BB\CC}

\title{Abstract integration, Combinatorics of Trees\\ and Differential Equations}
\author{Massimiliano Gubinelli\\[.4em]
{\small CEREMADE, Universit\'e Paris Dauphine, France}\\[-.2em]{\small  \texttt{gubinelli@ceremade.dauphine.fr}}}
\date{September 2008}

\begin{document}

\maketitle
\begin{abstract}
This is a review paper on recent work about the connections between rough path theory, the Connes-Kreimer Hopf algebra on rooted trees and the analysis of finite and infinite dimensional differential equation. We try to explain and motivate the theory of rough paths introduced by T. Lyons in the context of differential equations in presence of irregular noises. We show how it is used  in an abstract algebraic approach to the definition of integrals  over paths which involves a cochain complex of finite increments. In the context of such abstract integration theories we outline a connection with the combinatorics of rooted trees.  As interesting examples where these ideas apply we present two infinite dimensional dynamical systems:  the Navier-Stokes equation and the Korteweg-de-Vries equation.
\\[.5cm]

\end{abstract}

\section{Introduction}
\label{sec:intro}

Rough path theory has been developed by T. Lyons for the analysis
 of the map $\Phi : x \mapsto y$ which sends a vector-valued driving
signal $x$ to the solution of the differential equation
\begin{equation}
  \label{eq:diff-eq}
dy_t = \sum_a f_a(y_t) dx^a_t, \qquad y_0 = \overline y 
\end{equation}
where $\overline y$ is the inital condition. This
equation has a well defined meaning as a non-autonomous ODE when $x_t$ is differentiable in
$t$. In applications however it is interesting to consider more
general driving signals, e.g. when taking into account random
perturbations of dynamical systems where  $x_t$ is
the Brownian motion which is almost surely
non-differentiable. Indeed the application to stochastic differential
equations (SDEs) has been the main motivation in the developments of
the theory. The standard approach is then to understand
eq.~(\ref{eq:diff-eq}) as the integral equation
\begin{equation}
  \label{eq:ito-sde}
y_t = \overline{y} + \int_0^t \sum_a f_a(y_s) dx^a_s , 
\end{equation}
provide a well-defined meaning for the \emph{stochastic integral} in
the r.h.s. and the proceed to solve the equation by standard
fixed-point methods. Usually the integral above can be an It\^o
integral or a Stratonovich one, but other more exotic choices are
possible in other stochastic contexts (e.g. the Skorohod integral, the
``normal ordered'' integral, etc\dots). Whenever we speak of a solution $y$ to the
SDE~(\ref{eq:diff-eq}) we mean suitable random function $y$ for which
 $f(y)$ can be integrated against $dx$ in an appropriate
sense.

Lyons' basic observation was that, in the case of a smooth control $x$, the solution of the
differential equation is a well-behaved function of the \emph{iterated
integrals} $X$ of $x$:
\begin{equation}
  \label{eq:integrals}
X^{\overline a}_{ts} = \int_s^t \int_s^{u_n} \cdots
\int_s^{u_2} dx^{a_1}_{u_1} \cdots dx^{a_n}_{u_n}  
\end{equation}
where we denote with $\overline a$ the multi-index
$(a_1,\dots,a_n)$ with $|\overline a|=n$ its
length. In some sense the iterated integrals encode the local
behaviour of the path $x$ well enough to faithfully recover its effect
on the solution $y$. In a system perspective, the iterated itegrals
provide a canonical set of coordinates for the analysis of non-linear
systems much like the standard Fourier coefficients are natural
coordinates for linear ones.

Like Fourier coefficients, iterated integrals enjoy nice relations
upon concatenation of paths: given three times $s<u<t$ we have the
celebrated Chen relations~\cite{MR0454968} between iterated integrals:
\begin{equation}
  \label{eq:chen}
X^{a_1 \cdots a_n}_{ts} = X^{a_1 \cdots
  a_n}_{tu}+X^{a_1 \cdots a_n}_{us}+\sum_{k=1}^{n-1}X^{a_1
  \cdots a_k}_{tu} X^{a_{k+1} \cdots a_n}_{us}  .
\end{equation}
These non-linear equations plays a fundamental role in the development of the theory.

\bigskip
The plan of the paper is the following: in Sect.~\ref{sec:ren} we motivate the simplest instance of a rough path in the context of integration against smooth approximations of an irregular path. Sect.~\ref{sec:algint} we describe an algebraic approach to integration and its use to define integrals against a rough path. This approach is not standard from the point of view of rough paths. More conventional expositions of the theory can be found in~\cite{MR1654527,MR2036784} and~\cite{MR2053040}.  In Sect.~\ref{sec:trees} we describe the combinatorics of iterated integrals of very general type via the Connes-Kreimer Hopf algebra on rooted trees following~\cite{ramif}. This will ultimately allow to define and solve differential equation associated to such integrals. Finally in Sect.~\ref{sec:inf-dim} we apply the objects and the related combinatorics of algebraic integrals to the study of infinite dimensional differential equations via two prototypical examples: the 3d Navier-Stokes equation as studied in~\cite{MR2227041} and the 1d periodic Korteweg-de-Vries equation~\cite{kdv}.

\section{Renormalizable theories of integration} 
\label{sec:ren}

The analysis of the integral equation~(\ref{eq:ito-sde}) in the case
when $x$ is a non-differentiable function can be split
in two parts: the definition of the integral on the r.h.s. and the
fix-point argument. Let us concentrate on the first part. Assume given
a continuous function $x: [0,T] \to \RR^n$ which we assume only
$\gamma$-H\"older continuous, i.e. for which the following estimate holds
$$
|x_t-x_s| \le C |t-s|^\gamma, \qquad t,s \in [0,T]
$$
(where the least constant $C$ define the H\"older norm $\|x\|_\gamma$)
and for simplicity let us restrict in this section to the case $1/3 < \gamma < 1$.
Let us pose the problem of giving a ``reasonable'' definition of
the integral
$$
I[\varphi]_{ts} = \int_s^t \varphi(x_u) dx_u 
$$
for some smooth one-form $\varphi : \RR^n \to \RR^n$. Proceeding
by approximation we consider a family $\{x(\eps)\}_{\eps > 0}$ of
smooth approximations to the path $x$ and let
$$
I[\varphi](\eps)_{ts} = \int_s^t \varphi(x(\eps)_u) dx(\eps)_u. 
$$
It is easy to convince ourselves that in general we do not have any
control of these integrals if $x(\eps) \to x$ when $\eps \to 0$ in the
$\gamma$-H\"older norm. It is then a remarkable fact that all the
possible integrals obtained varying the function $\varphi$ all
converge at the same time (provided $\varphi$ is sufficiently smooth)
when the approximated iterated integral of order two:
$$
X(\eps)^{a_1 a_2}_{ts} = \int_s^t \int_s^u
dx(\eps)^{a_1}_v dx(\eps)^{a_2}_u
$$
converges to a function $X : [0,T]^2 \to
\RR^n \otimes \RR^n$ in the sense that
$$
\sup_{0 \le s < t \le T}\frac{|X(\eps)^{a_1
    a_2}_{ts}-X^{a_1 a_2}_{ts}|}{|t-s|^{2\gamma}} 
\to 0\qquad \qquad \text{as $\eps \to 0$.}
$$
So in some sense we can claim that the integration theory agains the
path $x$ is well defined as long as we are able to control the
convergence of the iterated integral of order two $X^{a_1 a_2}$: all the other
integrals (and more fundamentally also all the higher order iterated
integrals) will turn out to be nice functions of the data given by the
path $x$ and $X^{a_1 a_2}$. 

It is suggestive to understand this phenomenon as a very simple
example of ``renormalizable theory'' where all the quantities of
interests have well defined meanings as functionals of a finite number of
fundamental objects whose intrinsic determination remains outside the scope of the
theory itself.

In this limited context we face the appeareance of the simplest
non-trivial example of a \emph{rough path}: the couple $(X^a,X^{a_1
  a_2})$ where $X^a_{ts} = x^a_t - x^a_s$ is a
$\gamma$-rough path, i.e. a path and some additional information in
the form of ``iterated integrals'', for which a complete theory
of integration and differential equations can be constructed (as we
will see shortly). The fact that it is enough to
consider only the second order integral is due to our hypothesis that
$\gamma > 1/3$.

Note that once the limit has been taken the object $X^{a_1
  a_2}$ is no more an iterated integral (in a classical
sense) and can be characterized more abstractly by the following two
properties
\begin{enumerate}
\item
 the Chen relation:
 \begin{equation}
   \label{eq:chen-step-two}
X^{a_1
  a_2}_{ts} = X^{a_1
  a_2}_{tu} + X^{a_1
  a_2}_{us} + X^{a_1}_{tu} X^{ a_2}_{us}   
 \end{equation}
\item a regularity condition
 \begin{equation}
   \label{eq:reg-step-two}
\sup_{0 \le s < t \le T}\frac{|X^{a_1
    a_2}_{ts}|}{|t-s|^{2\gamma}}  < \infty.
 \end{equation}
\end{enumerate}
and more interestingly there could be  more than one
possible choice for this object compatible with these two conditions  leading to different integration
theories.

\section{Algebraic integration}
\label{sec:algint}

To understand how the iterated integrals $X$ comes into play in the
definition of $I[\varphi]$ we need some tool which allows us to
analyze the ``local'' (with respect to the parameter) behaviour of the
integral. Specifically we want to expand the integral in a short
interval as
\begin{equation}
  \label{eq:expansion-of-I}
I[\varphi]_{ts} = \varphi_a(x_s) X^a_{ts} + \partial_b
\varphi_a(x_s) X^{b a}_{ts} + r_{ts}  
\end{equation}
where $r$ stands for some remainder term which we hope will be of
higher order in $|t-s|$ than the other terms. The main property of the
integral $I[\varphi]$ is its trivial behaviour under splitting of the
integration interval: $I[\varphi]_{ts} = I[\varphi]_{tu} +
I[\varphi]_{us}$ for $s < u < t$, this of course means also that the
integral can be written as the increment of a function
$I[\varphi]_{ts} = f_t - f_s$ (for example taking $f_t =
I[\varphi]_{0t}$). In~\cite{MR2091358} we introduced a cochain complex
$(\CC_*,\delta)$ which encodes this basic property of integrals. For
$n \ge 1$ let
$\CC_n$ be set of functions $a\in C([0,T]^n;\RR)$  such that $a_{t_1
  \cdots t_n} = 0$ if $t_i = t_{i+1}$ for some $1 \le i \le
n-1$. Elements in $\CC_n$ will be called \emph{$n$-increments}. 
Define a coboundary $\der : \CC_n \to \CC_{n+1}$ as $\der a_{ts}
= a_t-a_s$ for $a \in \CC_1$, $\der a _{tus} = a_{ts}-a_{tu}-a_{us}$
for $a\in \CC_2$ and so on. Let $\ZCC_n = \textrm{Ker} \der \cap
\CC_n$  and $\BCC_n = \textrm{Im} \der \cap \CC_n$. It is easy then to verify that $\der \der
= 0$ and that the complex $(\CC_*,\der)$ is exact, i.e. that
$\text{Im} \der = \text{Ker} \der$ at any $\CC_n$, $n \ge 1$. In
particular a $2$-increment $a\in\CC_2$ is  the increment of a function
$f \in \CC_1$ if and only if $\der a = 0$. So at $\CC_2$ the
coboundary measures the degree of "exactness'' of 2-increments.
Moreover a key fact is the following: if $\der a$ is suitably small
then there exists \emph{only one} function $f$ (modulo constants) such
that
\begin{equation}
  \label{eq:basic-decomp}
\der f = a + r
\end{equation}
where the remainder is small. To be more precise we need to introduce
the relevant notion of ``smallness''. We say that $a\in \CC_2^\gamma$
if 
$$
\|a\|_\gamma := \sup_{s<t} \frac{|a_{ts}|}{|t-s|^\gamma} < \infty .
$$
And similarly $b \in \CC_3^\gamma$ if
$$
\|b\|_\gamma := \sup_{s<u<t} \frac{|b_{tus}|}{|t-s|^\gamma} < \infty .
$$
Both $\CC_2^{\gamma}$ and $\CC_3^\gamma$ are Banach spaces when
endowed with the norms $\|\cdot\|_\gamma$. Define $\CC_n^{1+} =
\cup_{z > 1} \CC_n^z$.  Moreover if $g \in \CC_n$ and $h \in \CC_m$ then we write $gh$ for the element of $\CC_{n+m-1}$ such that
$(gh)_{t_1 \cdots t_{n+m-1}} = g_{t_1 \cdots t_n} h_{t_n \cdots t_{n+m-1}}$.

\bigskip

The key result in this theory of increments is the existence of a \emph{sewing map} which provide a natural inverse operation to $\der$:

\begin{theorem}[The sewing map] There exists a unique bounded and linear
  map $\Lambda : \ZCC_3^{1+} \to \CC^{1+}_2$ such that $\der \Lambda h = h$ for any $h \in \ZCC_3^{1+}$.
\end{theorem}

In particular the map $(1-\Lambda \der) : \CC_2 \to \CC_2 $ projects
(in essentially a unique way)
suitable elements of $\CC_2$ to $\ZCC_2 = \BCC_1$: the space of
increments of functions.
Then going back to the decomposition~(\ref{eq:basic-decomp}) and
assuming that  $\der a\in \CC_3^{1+} $ we can form the exact 2-increment
$a - \Lambda \der a = \der f$ and obtain that the remainder $r$ is
given by $\Lambda \der a \in \CC_2^{1+}$. It is also easy to
see that, the decomposition~(\ref{eq:basic-decomp}) is unique if we
require $r \in \CC_2^{1+}$, indeed if two such decomposition
exists, $(f,r)$ and $(f',r')$, their difference satisfy $\der(f-f')=r-r'
\in \CC_2^z$ for some $z > 1$ and $f-f'$ would be a function of
$z$-H\"older class ans since $z>1$ that this function would take the
constant value zero.

With this notions at hand we realize that eq.~(\ref{eq:expansion-of-I}) is nothing more
than an instance of a decomposition similar
to~(\ref{eq:basic-decomp}). As we have already seen we can determine
both $I[\varphi]$ and $r$ at once using only the well-defined
remaining terms in the
r.h.s., so letting $a_{ts} =  \varphi_a(x_s) X^a_{ts} + \partial_b
\varphi_a(x_s) X^{b a}_{ts} $, in order to be able to
apply $\Lambda$, we need to require that
$\der a \in \CC_3^{1+}$. By some easy computation using
the fact that $\der$ satisfy some kind of Leibniz rule, we
get
$$
\der a_{tus} =  - \der \varphi_a(x)_{tu} X^a_{us} -  \der \partial_b
\varphi_a(x)_{tu} X^{b a}_{ts} +  \partial_b
\varphi_a(x_s) \der X^{b a}_{tus}
$$
If we now exploit the Chen relation for $X^{b a}$ we can
simplify this expression further
\begin{equation}
\label{eq:step-2-obstruction}
\der a_{tus} =  - [\der \varphi_a(x)_{tu} -  \partial_b
\varphi_a(x_s) \der X^{b}_{tu} ]X^a_{us} -  \der \partial_b
\varphi_a(x)_{tu} X^{b a}_{ts} 
\end{equation}
Assuming that $\varphi$ is sufficiently smooth ($C^2$ is enough) we
obtain easily that
$$
\der \varphi_a(x) -  \partial_b
\varphi_a(x) \der X^{b} \in \CC_2^{2\gamma}, \qquad
 \der \partial_b
\varphi_a(x) \in \CC^\gamma_2
$$
so that taking into account all the regularities we end up with $\der
a \in \CC^{3\gamma}_3$. Then our assumption on $\gamma$ ensure that
$3\gamma > 1$ and we can prove that there exists a unique couple
$(f,r)$ with $f \in \CC_1$ and $r \in \CC_2^{1+}$ such
that
\begin{equation}
\label{eq:renorm-int}
\der f = \varphi_a(x_s) X^a_{ts} + \partial_b
\varphi_a(x_s) X^{b a}_{ts} + r_{ts}.
\end{equation}
By construction $f$ depends only on the $\gamma$-rough path $X$ and on
$\varphi$. We can then  \emph{define} the integral by $$I[\varphi] =
\der f=   (1-\Lambda \der)[ X^a \varphi_a(x)  + X^{b a} \partial_b
\varphi_a(x) ].$$

To motivate the fact that the decomposition~(\ref{eq:renorm-int}) is a sort of renormalized integral we can make the following observation. Take a partition $\{t_i\}$ of $[s,t]$ of size $\Delta$ and consider the Riemman sums $S_\Delta = \sum_i \varphi(x_{t_i}) X_{t_{i+1} t_i}$. In general we have no mean to say that these sums are convergent as $\Delta \to 0$. However given a $\gamma$-rough path $X$ we can perform a subtraction to these sums and define
$$
S'_\Delta = \sum_i \left[ \varphi_a(x_{t_i}) X^a_{t_{i+1} t_i} +  \partial_b
\varphi_a(x_{t_i}) X^{b a}_{t_{i+1} t_i}\right]
$$
then using the decomposition\ref{eq:renorm-int} as $\Delta \to 0$ we have the limit
$$
S'_\Delta = \sum_i (\der f)_{t_{i+1} t_i} - \sum_i r_{t_{i+1} t_i} = (\der f)_{t s} - \sum_i o(|t_{i+1} t_i|) \to (\der f)_{t s} 
$$
since the first sum telescopes and the second is easily show to converge to zero. 

Moreover the regular dependence on the data ensure that smooth approximations
 $I[\varphi](\eps)$ converge to $I[\varphi]$ as here defined as long
 as we can prove the convergence of approximating path (and its second
 order iterated integrals) to the rough path $X$.

Exploiting the sewing map, the integral equation~(\ref{eq:ito-sde}) in presence of a
$\gamma$-rough path with $\gamma > 1/3$ can be
understood as a fixed-point equation for an unknown $y \in \CC_1$:
$$
\der y = (1-\Lambda \der)[ X^a \varphi_a(y) + 
 X^{a c} \partial_b \varphi_a(y) \varphi^b_{c}(y) ]
$$
these equation can be solved by a standard iteration method in a
suitable subspace of $\CC_1$ (\cite{MR2091358} for details).

It is also possible to construct higher order
iterated integrals starting from low order ones where for iterated
integrals we means simply object which obey Chen's
relations~(\ref{eq:chen}). Take for example the third order object
$X^{a_1 a_2 a_3} \in \CC_2$, eq.~(\ref{eq:chen}) can be
written as a statement about the coboundary of $X^{a_1 a_2 a_3}$:
$$
\der X^{a_1 a_2 a_3}_{tus} = X^{a_1 a_2
  }_{tu} X^{a_3}_{us} + X^{a_1}_{tu} X^{ a_2 a_3}_{us} 
$$
and in the above hypothesis on $X$ it is easy to check that the
r.h.s. belongs to $\ZCC_2^{3\gamma} \subset \ZZ\CC_2^{1+}$ so that it is in the domain of
$\Lambda$ and we can \emph{define}
$$
X^{a_1 a_2 a_3} := \Lambda[X^{a_1 a_2
  } X^{a_3} + X^{a_1} X^{ a_2 a_3}] 
$$
as the unique solution in $\CC_2^{1+}$ of the Chen
relation. Iteratively this allows to construct all the higher order
iterated integrals. Some more work allows to prove a  uniform estimate
on the growth of the norms involved in the procedure~\cite{MR1654527,MR2091358}:
\begin{equation}
\label{eq:growth-lin}
\|X^{\overline a}\|_{|\overline a|\gamma} \le C_1
\frac{C_2^{|\overline a|}}{(|\overline a|!)^{\gamma}}
\end{equation}
for any multiindex $a$. 

\section{Trees}
\label{sec:trees}

What happens if $\gamma \le 1/3$? The obstuction to the exactness of
the increment $a$ in eq.~(\ref{eq:step-2-obstruction}) will no more
belong to the domain of the sewing map. Indeed in
eq.~(\ref{eq:expansion-of-I}) we cannot anymore expect that the
remainder belongs to $\CC^{1+}_2$ and at the very least this would
affect our argument for uniqueness. We are then forced to proceed and
expand further the integral, or from another point of view, to add
some counterterms to remove large contributions to $r$. 

To understand the general structures that we need if we proceed
further in the expansions it is better to take as working bench the
more difficult case of the differential equation
(\ref{eq:ito-sde}). The series solution will be indexed by rooted
trees: a phenomenon which is present already for solutions to the ODE
$dy/dt = f(y)$~\cite{MR0305608,MR2106008}.

\subsection{Labelled rooted trees}
Given a finite set $\cl$, the $\cl$-labeled rooted trees are the
finite graphs where labels of $\cl$ are attached to each vertex and
where there is a special vertex called \emph{root} such
that there is a unique path from the root to any other vertex. 

Some examples of rooted trees labeled by $\cl=\{1,2,3\}$ are 
$$
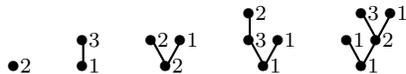
\tnodel 2 \qquad 
\pstree{\tnodel 1}{\tnodel 3}
\qquad
\pstree{\tnodel 2}{\tnodel 2 \tnodel 1}
\qquad
\pstree{\tnodel 1}{\pstree{\tnodel 3}{\tnodel 2} \tnodel 1}
\qquad
\pstree{\tnodel 1}{\tnodel 1 \pstree{\tnodel 2}{\tnodel 3 \tnodel 1}} 
$$
 Trees
does not distinguish the order of branches to each vertex.
Given $k$ $\cl$-decorated rooted trees $\tau_1,\cdots,\tau_k$ and a label $a \in\cl$ we define
$\tau = [\tau_1,\cdots,\tau_k]_a$ as the tree obtained by attaching the $k$
roots of $\tau_1,\cdots,\tau_k$ to a new vertex with label $a$ which will be the root
of $\tau$. Any decorated rooted tree can be constructed using the
simple decorated tree $\bullet_a$ ($a \in\LL$)
and the operation $[\cdots]$, e.g.
$$
[\bullet] = \pstree{\tnode}{ \tnode}
\qquad
[\bullet,[\bullet]] = \pstree{\tnode}{\pstree{\tnode}{\tnode} \tnode},
\qquad \text{etc\dots}
$$

Denote $\TT_\cl$ the set of all $\cl$
decorated rooted trees and let $\TT$ the set of rooted trees without
decoration (i.e. for which the set of labels $\cl$ is made of a single
element). There is a canonical map $\TT_\cl \to \TT$ which simply
forget all the labels and every function on $\TT$ can be extended,
using this map to a function on $\TT_\cl$ for any set of labels
$\cl$. 
Let $|\cdot| : \TT \to \RR$ the weight which counts the
number of vertices of the (undecorated) tree defined recursively as
$$
|\bullet | = 1, \qquad |[\tau_1,\dots,\tau_k]| = 1+|\tau_1|+\cdots+|\tau_k|
$$
moreover the \emph{factorial} $\gamma(\tau)$ of a tree $\tau$ is
$$
\gamma(\bullet) = 1, \qquad \gamma([\tau_1,\dots,\tau_k]) =
|[\tau_1,\dots,\tau_k]| \gamma(\tau_1) \cdots \gamma(\tau_k)
$$
finally the \emph{symmetry factor} $\sigma$ is given by the
recursive formula $\sigma(\tau) = 1$ for $|\tau| = 1$ and
\begin{equation}
  \label{eq:sigma-prop}
\sigma([\tau^1 \cdots \tau^k]_a) = \frac{k!}{\delta(\tau^1,\dots,\tau^k)} \sigma(\tau^1) \cdots \sigma(\tau^k)  
\end{equation}
where $\delta(\tau^1,\cdots,\tau^{k})$ counts the number of different ordered $k$-uples $(\tau^1,\cdots,\tau^{k})$ which corresponds to the same (unordered) collection  $\{\tau^1,\cdots,\tau^{k}\}$ of subtrees.
The factor $k!/\delta(\tau^1,\dots,\tau^k)$ counts the order of the
subgroup of permutations of $k$ elements which does not change the
ordered $k$-uple $(\tau^1,\cdots,\tau^{k})$. Then $\sigma(\tau)$ is
 is the order of the subgroup of permutations on the vertex of
the tree $\tau$ which do not change the tree (taking into account also
the labels). Another equivalent recursive definition for $\sigma$ is
$$
\sigma([(\tau^1)^{n_1}\cdots(\tau^k)^{n_k}]_a) = n_1!\cdots n_k!
\sigma(\tau^1)^{n_1}\cdots \sigma(\tau_k)^{n_k}
$$
where $\tau^1,\dots,\tau^k$ are distinct subtrees and $n_1,\dots,n_k$
the respective multiplicities.

The algebra $\cA \TT_\cl$ is the commutative free polynomial algebra
generated by $\{1\}\cup \TT_\cl$ over $\RR$, i.e. elements of $\cA \TT_\cl$ are
finite linear combination with coefficients in $\RR$ of formal monomials in the form $\tau_1
\tau_2 \cdots \tau_n$ with $\tau_1,\dots,\tau_n \in \TT_\cl$ or of the
unit $1\in \cA\TT_\cl$. The tree monomials are called \emph{forests} and
are collectively denoted $\FF_\cl$; we include the empty forest $1 \in
\FF_\LL$.
 The algebra $\cA \TT_\cl$
 is endowed with a graduation $g$ given  by $g(\tau_1
\cdots \tau_n) = |\tau_1|+\cdots+|\tau_n|$ and $g(1) = 0$. This graduation induces a
corresponding filtration of $\cA \TT_\cl$ in finite dimensional linear
subspaces $\cA_n \TT_\cl$ generated by the set $\FF_\cl^n$ of forests of degree $\le n$.

Any map $f : \TT_\cl \to
A$ where $A$ is some commutative algebra, can be extended in a unique way to
a homomorphism $f : \cA
\TT_\cl \to A$ by setting:
$
f(\tau_1 \cdots \tau_n) = f(\tau_1) f(\tau_2) \cdots  f(\tau_n)
$.

In the following we will use letters $\tau,\rho,\sigma,\dots$ to
denote trees in $\TT_\cl$ or forests  in $\FF_\LL$, the degree $g(\tau)$ of a forest $\tau \in \FF_\cl$ will also be written as $|\tau|$. Roman letters $a,b,c,\dots \in \LL$ will denote vector  indexes (i.e. labels) while $\ol{a},\ol{b},\dots$ will denote multi-indexes with values in $\LL$: $\overline{a}=(a_1,\dots,a_n) \in \LL^n$ with $|\overline a|=n$ the size of this multi-index.

\subsection{Iterated integrals and the Connes-Kreimer Hopf algebra}
Given a smooth path $x \in C^1([0,T],\RR^n)$ we can canonically
associate to it a family of 2-increments $X^\tau$ indexed by trees
labelled by $\LL = \{1,\dots,n\}$ by the iterated integrals
\begin{equation}
  \label{eq:tree-integrals}
X^{\troot_a}_{ts} = x_t-x_s, \qquad X^{[\tau_1 \cdots \tau_n]_a}_{ts}
= \int_s^t X^{\tau_1}_{us} \cdots  X^{\tau_n}_{us} dx^a_u
\end{equation}
Iterated integrals like~(\ref{eq:integrals}) corresponds to ``linear''
trees $\tau = [[\cdots [\troot_{a_n}]_{a_{n-1}}
\cdots]_{a_2}]_{a_1}$ and by abuse of notation we will
continue to write $X^{a_1 \cdots a_n}$ for such $X^\tau$.
The generalization of the Chen's relations~(\ref{eq:chen}) involves
the Hopf algebra structure on on labelled
trees essentially introduced by Connes and Kreimer which we now
brefly describe.

On the algebra $\cA \TT_\cl$ we can define a counit $\eps: \cA
\TT_\cl \to \RR$ as an algebra homomorphism such that $\eps(1) = 1$
and $\eps(\tau) = 0$ otherwise and a \emph{coproduct} $\Delta :
\cA\TT_\cl \to \cA\TT_\cl \otimes \cA\TT_\cl$ as the algebra
homomorphism
such that
\begin{equation}
  \label{eq:delta-def}
\Delta(\tau) = 1\otimes \tau  + \sum_{a \in \cl} (B^a_+ \otimes \text{id} )[\Delta(B^a_-(\tau))] 
\end{equation}
on trees $\tau \in \TT_\LL$,
where $B_+^a(1) = \troot_a$ and $B_+^a(\tau_1 \cdots \tau_n) = [\tau_1 \cdots \tau_n]_a$ and
$B_-^a$ is the inverse of $B_+^a$ or is equal to zero if the tree
root does not have label $a$, i.e. 
\begin{equation*}
B_-^a(B^b_+(\tau_1 \cdots \tau_n))= \begin{cases}
 \tau_1 \cdots \tau_n & \text{if $a=b$}\\
0 & \text{otherwise.}
\end{cases}
\end{equation*}
%
Endowed with $\eps$ and $\Delta$ the algebra $\cA \TT_\cl$ become a
bialgebra, there exists also an antipode $S$ which complete the
definition of the Hopf algebra structure on $\TT_\cl$ as described by
Connes-Kreimer~\cite{MR1660199} (in the unlabeled case). 

We will often use Sweedler's notation for the coproduct $\Delta \tau = \sum \tau_{(1)} \otimes \tau_{(2)}$ but we also introduce a counting function $c : \TT_\cl \times \TT_\cl \times \FF_\cl \to \NN$ such that
$
\Delta \tau = \sum_{\rho \in \TT_\cl, \sigma \in \FF_\cl} c(\tau,\rho, \sigma) \rho \otimes \sigma
$
moreover it will be useful to consider also the reduced coproduct
$\Delta' \tau = \Delta \tau - 1 \otimes \tau - \tau \otimes 1$ with
counting function $c'$.

The generalization of eq.~(\ref{eq:chen}) reads as follows:

\begin{theorem}[Tree multiplicative property]
\label{th:equiv}
The map $X$ satisfy the algebraic relation
\begin{equation}
  \label{eq:alg-relations}
\der X^{\sigma} = X^{\Delta'(\sigma)} , \qquad \sigma \in \cA\TT_\cl
\end{equation}
\end{theorem}

Let us give an example in one dimension ($d=1$) where trees are not
decorated. The
forests with $|\tau| \le 3$ are
\smalltrees
\begin{equation*}
\tsnode, \aabb, \tsnode \tsnode, \aaabbb, \tsnode \aabb, \tsnode \tsnode \tsnode, \aababb  
\end{equation*}
The action of the reduced coproduct on these forests and the
corresponding action of the coboundary on the iterated integrals are
given by

\begin{minipage}{.45\textwidth}
\begin{equation*}
\Delta'{ \aabb } = {\tsroot} \otimes { \tsroot }
\end{equation*}
\begin{equation*}
\Delta' ( {\tsroot \tsroot} ) = 2 {\tsroot} \otimes {\tsroot  }
\end{equation*}
\begin{equation*}
\Delta'   \aaabbb = \aabb \otimes \tsroot + \tsroot \otimes \aabb 
\end{equation*}
\begin{equation*}
\Delta' ({ \tsroot \aabb} ) = \tsroot \otimes \tsroot \tsroot + \tsroot \tsroot \otimes \tsroot
+ \aabb \otimes \tsroot + \tsroot \otimes \aabb
\end{equation*}
\begin{equation*}
\Delta'(\tsroot^3) = 3 \tsroot^2 \otimes \tsroot + 3 \tsroot \otimes \tsroot^2  
\end{equation*}
\begin{equation*}
 \Delta' \aababb  = \tsroot \otimes \tsroot \tsroot + 2 \aabb \otimes \tsroot
\end{equation*}
\end{minipage}
\begin{minipage}{.45\textwidth}
\begin{equation*}
\der X_{tus}^{ \aabb } = X_{tu}^{\tsroot} X_{us}^{ \tsroot }
\end{equation*}
\begin{equation*}
\der X_{tus}^{\tsroot \tsroot}  = 2 X_{tu}^{\tsroot} X_{us}^ {\tsroot  }
\end{equation*}
\begin{equation*}
\der X_{tus}^{ \aaabbb } = X_{tu}^{\aabb}X_{us}^{\tsroot} + X_{tu}^{\tsroot} X_{us}^{\aabb} 
\end{equation*}
\begin{equation*}
\der X_{tus}^{ \tsroot \aabb}  = X_{tu}^{\tsroot} X_{us}^{ \tsroot \tsroot} + X_{tu}^{\tsroot \tsroot} X_{us}^{ \tsroot}
+ X_{tu}^{\aabb} X_{us}^{\tsroot} + X_{tu}^{\tsroot}X_{us}^{ \aabb}
\end{equation*}
\begin{equation*}
\der X_{tus}^{\tsroot^3} = 3 X_{tu}^{\tsroot^2} X_{us}^{ \tsroot } + 3 X_{tu}^{\tsroot}X_{us}^{\tsroot^2}  
\end{equation*}
\begin{equation*}
 \der X_{tus}^{\aababb}  = X_{tu}^{\tsroot}X_{us}^{\tsroot \tsroot} + 2 X_{tu}^{\aabb}X_{us}^{ \tsroot}
\end{equation*}
\end{minipage}
\medskip

As a further example
consider the iterated integrals $T^\tau$ associated to the identity path 
$$
T^{\troot}_{ts} = t-s, \qquad T^{[\tau_1 \cdots \tau_n]}_{ts} = \int_s^t T^{\tau_1}_{us}\cdots T^{\tau_n}_{us} du
$$  
By induction it is not difficult to prove that
$
T^{\tau}_{ts} = (t-s)^{|\tau|} (\tau!)^{-1}
$,  so applying Thm.~\ref{th:equiv} to $T^\tau$ we get a remarkable
binomial-like formula for the Connes-Kreimer coproduct
\begin{equation}
  \label{eq:tree-binomial}
  (a+b)^{|\tau|} = \sum  \frac{\tau!}{\tau^{(1)}! \tau^{(2)}!}
  a^{|\tau^{(1)}|} b^{|\tau^{(2)}|} .
\end{equation}

\bigskip
Iterated integrals of this sort appears naturally when trying to expand in series the solution of driven differential equations:
for any analytic vectorfield $f : \RR^d \to \RR^n$ and any smooth path
$x \in C^1([0,T],\RR^n)$, the solution of the differential equation
$
 dy_t = \sum_{a\in\LL}  f_a(y_t) dx_t^a$, with $y_0 = \eta
$  
admit locally the series representation
\begin{equation}
\label{eq:tree-rep}
\der y_{ts} = \sum_{\tau \in \TT_\LL} \frac{1}{\sigma(\tau)} \phi^{f}(\tau)(y_s) X^\tau_{ts}, \qquad y_0 = \eta
\end{equation}
where   we recursively define
functions $\phi^f(\tau)$ as
\begin{equation*}
\phi^f(\troot_a)(\xi) = f_a(\xi), \qquad
\phi^f([\tau^1 \cdots \tau^k]_a)(\xi) =  \sum_{\ol b  \in \II \LL_1: |\ol b| = k}
f_{a;b_1\dots b_k}(\xi) \prod_{i=1}^{k} [\phi^{f}(\tau^i)(\xi)]^{b_i}  .
\end{equation*}
%

\bigskip

Note that using the regularity of the path $x$ the iterated integrals
$X^{\tau}$ can be simplified: indeed Chen~\cite{MR0454968} proved that
products of iterated integrals can be always expressed as linear
combination of iterated integrals via the \emph{shuffle product}:
\begin{equation}
  \label{eq:shuffle}
X^{a_1 \cdots a_n }_{ts}  X^{b_1 \cdots b_m}_{ts} =
\sum_{\overline c \in \text{Sh}(\overline a, \overline b)}
X^{c_1\cdots c_{n+m}}_{ts}  
\end{equation}
where given two multi-indexes $\overline a = (a_1,\cdots, a_n)$ and
$\overline b = (b_1, \cdots, b_n)$ their \emph{shuffles} $\text{Sh}(\overline a, \overline b)$ is the set of all the possible permutations of the $(n+m)$-uple $(a_1,\dots,a_n,b_1,\dots,b_m)$ which does not change the ordering of the two subsets $\overline a$, $\overline b$.
Using relation~(\ref{eq:shuffle}) every $X^{\tau}$ can be reduced to a linear combination of standard iterated
integrals $\{ X^{\ol a}
\}_{\ol a}$. 

\bigskip
It is however interesting that many costructions related to integrals
and solution of driven differential equations do not depend on
eq.~(\ref{eq:shuffle}) being valid. For example in the theory of the
It\^o integral eq.~(\ref{eq:shuffle}) does not hold and more complex
algebraic relations have to be considered: if $x$ is a multidimensional
Brownian motion and $X$ defined via It\^o integration we have
$$
X^a_{ts} X^b_{ts} = X^{ab}_{ts}+X^{ba}_{ts} + \delta_{ab} (t-s).
$$

\bigskip

Let us clarify the algebraic framework in which we consider possible
integration theories.
It turns out that the only data we really need to build a family
$\{X^\tau\}_{\tau \in \TT_\LL}$ satisfying~(\ref{eq:alg-relations}) is
given by a set of linear maps $\{I^a\}_{a \in
  \LL}$ defined on a domain $\DD_I \subset \CC^+_2$ (the unital
algebra obtained extending $\CC_2$ with a unit $e$ and considering the point-wise
product) to $\CC_2$ satisfying two properties:
\begin{enumerate}
\item 
$
I(hf)_{ts} =  I(h)_{ts} f_s$ for all $ h \in \DD_I, f \in \CC_1
$ and
where $(hf)_{ts} = h_{ts} f_s$;
\item
$
\der I(h)_{tus} = I(e)_{tu}h_{us} + \sum_i I(h^{1,i})_{tu} h^{2,i}_{us}$ when $h\in
  \DD_I$ with $ \der h_{tus} = \sum_i h^{1,i}_{tu} h^{2,i}_{us}$ and $h^{1,i} \in \DD_I$ .
\end{enumerate}
 
Then
given a family $\{I^a\}_{a\in\LL}$ of such integral maps on a common
algebra $\DD\subseteq \CC_{2}$ we can associate to them a family $\{X^{\tau}\}_{\tau\in\FF_\LL}$ recursively as
$$
X^{\troot_a}_{ts} = I^a(e)_{ts}, \qquad X^{[\tau^1 \cdots \tau^k]_a}_{ts} = I^a(X^{\tau^1 \cdots \tau^k})_{ts}, \qquad X^{\tau^1\cdots \tau^k}_{ts}= X^{\tau^1}_{ts} \cdots  X^{\tau^k}_{ts}.
$$   
In this way we estabilish an algebra homomorphism from $\cA\TT_\LL$ to
a subalgebra of $\CC_2$ generated by the $X^\tau$-s. This homomorphism
send the operation $B_+^a$ on $\cA\TT_\LL$ to the integral map $I^a$
on $\CC_2$ in such a way that Theorem~\ref{th:equiv} holds.

Given $X^\tau$ for any $\tau\in\TT_\LL$ with $|\tau| \le n$ and such
that $\|X^\tau\|_{\gamma |\tau|} \le C$ for some constant $C$ and
$\gamma \in (1/(n+1),1)$ we are able to extend to any $\tau$ with
$|\tau| > n$ the ``integrals'' $X^\tau
\in \CC^{\gamma |\tau|}_2$ in a unique manner solving  the
equation $\der X^{\tau} = X^{\Delta ' \tau}$ with the aid of the
sewing map. This extension will satisfy the bound
\begin{equation}
  \label{eq:branched-bounds}
\|X^\tau\|_{\gamma |\tau|} \le C_1 C_2^{|\tau|} q_\gamma(\tau)  
\end{equation}
for any $\tau$ with $C_1,C_2$ two finite constants and
$q_\gamma(\tau)$ a function satisfying
$$
q_\gamma(\tau) = \frac{1}{2^{\gamma |\tau|}-2} \sum'
q_\gamma(\tau^{(1)}) q_\gamma(\tau^{(2)}).
$$
The actual asympotic behavior of $q_\gamma$ for $|\tau| \to \infty$ is
not yet know but we conjecture that it should hold the equivalence
\begin{equation}
  \label{eq:q-conj}
  q_{\gamma}(\tau) \simeq (\tau !)^{-\gamma}
\end{equation}
(see~\cite{ramif}).
This asymptotic behavior is satisfied on the subset
of linear trees and is consistent with our
results on  tree-indexed iterated integrals in the context of 3d
Navier-Stokes equation~\cite{MR2227041} (see also Sect.~\ref{sec:inf-dim}).

We call the object  $\{X^{\tau}\}_{\tau}$ satisfying~(\ref{eq:alg-relations})
and~(\ref{eq:branched-bounds}) for some $\gamma \in (0,1)$ a
($\gamma$-)\emph{branched rough path}. By the above consideration it
is clear that it is determined by the finite subset
$\{X^{\tau}\}_{\gamma |\tau| \le 1}$.

Assuming that the integrals $\{I^a\}_{a \in \LL}$ generate a
$\gamma$-branched rough path in the above sense, then we can extend
the integral to a larger class of integrands. Take a set of constants
$\{h^\tau_0 \}_\tau$ and consider the path $h \in \CC_1$ defined as
$
h_t = \sum_\tau h^\tau_0 X^{\tau}_{t0}
$, then
\begin{equation}
\label{eq:controlled-h}
\der h_{ts} = \sum_\tau h^\tau_0  \sum' X^{\tau^{(1)}}_{ts}
X^{\tau^{(2)}}_{s0} = \sum_{\tau,\rho,\sigma} c'(\tau,\rho,\sigma) h^\tau_0 X^{\rho}_{ts}
X^{\sigma}_{s0} = \sum_\tau h^{\tau}_s X^{\tau}_{ts}
\end{equation}
where we introduced paths $h^\tau \in \CC_1$ by
$
h^\tau_s =  \sum_{\sigma,\rho,\tau} c'(\rho,\tau,\sigma) h^\rho_0 X^{\sigma}_{s0}
$.
We must have
$$
0= \der \sum_\tau  X^{\tau} h^{\tau} = \sum_\tau X^{\tau^{(1)}}
X^{\tau^{(2)}} h^\tau - X^{\tau} \der h^\tau 
$$
so we can check that 
\begin{equation}
  \label{eq:controlled-derivatives}
\der h^{\tau} =  \sum_{\rho,\tau,\sigma} c'(\rho,\tau,\sigma) 
X^{\sigma} h^\rho   
\end{equation}
also holds. 
At this level these relations are formal since they require an
infinite number of terms to hold, however exploiting our analysis of
the $(\CC_*,\der)$ complex and the sewing map we can work modulo
$\CC_2^{1+}$ and neglect the terms in the expansions which involve
$X^\tau$-s with $\gamma|\tau| > 1$. A more accurate analysis (which
can be found in~\cite{ramif}) reveals that for the purpose of
integrating and solving
differential equations we can actually work modulo larger terms but
for the sake of clarity we refrain to do this here.

We call a \emph{controlled path} (by $X$) any
 path $h$ which satisfy modulo $\CC^{1+}_2$ the eq.~(\ref{eq:controlled-h}) together with
eq.~(\ref{eq:controlled-derivatives}) for all trees $\tau$ with $\gamma |\tau| < 1$.  The
integral $I$ can then be extended to integrate any controlled path. Indeed
since $h e = e h + \der h$ it
is a consistent definition to set
\begin{equation*}
  \begin{split}
I^a(h) & = I^a(he) = I^a(e) h + I^a(\der h) = X^{\troot_a} h +
\sum_{\gamma |\tau|< 1} I^a(X^\tau) h^\tau + I^a(\crr)
\\ &  = X^{\troot_a}
h + \sum_{\gamma |\tau|< 1} X^{[\tau]_a} h^\tau +\crr
\end{split}
\end{equation*}
where we exploited the properties of the integral $I^a$ and the
definition of $X$ and where $\crr\in \CC^{1+}_2$ is a generic remainder
term (possibly different from line to line). Controlled paths are also
stable under mapping by sufficiently smooth functions $f$ (the degree
of differentiability depends on $\gamma$): i.e. if $y$ is a controlled
path, then $z=f(y)$ is again a controlled path with an explicit
formula for the coefficients $z^\tau$ in term of the derivatives of
$f$ and of the coefficients $y^\tau$. Using these properties we can
consider the differential equation for the integrals $I$ and
vectorfield $f$: 
$$
\der y = I^a(f_a(y)), \qquad y_0 = \eta
$$
in the space of controlled paths $y$ and solve it by a fixed-point argument.

\section{Infinite dimensional dynamical systems}
\label{sec:inf-dim}

In this section we would like to show how the ideas and the tools
described before could be applied in the context of infinite
dimensional dynamical system by introducing operator-valued 
iterated integrals (and rough paths) which by their non-commutative nature are intrinsically indexed by trees. In particular we will discuss two
different examples:  the 1d periodic
Korteweg-de-Vries (KdV) equation and
the 3d Navier Stokes (NS) equation. In the first case 
we will exploit the increment complex to
analyze perturbatively the solution for irregular initial data
, in the second case 
we will use series
expansion over trees of the solution to analyze the dynamics for large
times
.

\medskip
As before we will restrict ourself to an overview of the results and
to a sketch of the arguments. The
interested readers can find a rigorous discussion
elsewhere~\cite{kdv,MR2227041} and he can refer
to~\cite{TindelGubinelli} for extension of these consideration covering the analysis of stochastic partial differential
equations driven by irregular noises.

\subsection{The KdV equation}
The 1d periodic KdV equation is the partial differential equation
\begin{equation}
  \label{eq:kdv-real}
\partial_t u(t,\xi) + \partial^3_\xi u(t,\xi) +
\frac12 \partial_\xi u(t,\xi)^2 = 0,\quad u(0,\xi) = u_0(\xi), \qquad
(t,\xi)\in\RR\times \bTT  
\end{equation}
on the torus $\bTT = [-\pi,\pi]$. When the initial data is not smooth
this equation must be interpreted in the integral form
\begin{equation*}
u(t) = U(t) u_0 + \frac{1}{2}\int_0^t U(t-s) \partial [u(s)]^2 ds  
\end{equation*}
where $U$ is the ``free propagator'' given by the linear equation
$
\partial_t U(t) = \partial^3 U(t)
$, $\partial$ denoting the spatial derivative, and where we consider the solution $u$ as a path in a space of
functions.
By going to the interaction picture $\tilde v_t = U(-t) u(t)$ we get
\begin{equation*}
\tilde  v_t = u_0 + \frac{1}{2}\int_0^t U(-s) \partial [U(s)\tilde v_s]^2 ds  
\end{equation*}
Then the Fourier coefficients $\{ v_t(k) \}_{k \in \ZZ}$  of $\tilde
v_t$ satisfy the equation
\begin{equation}
  \label{eq:kdv-base}
 v_t(k) =   v_0(k) + \frac{ik}2 \sum_{k=k_1+k_2, k_n \neq 0} \int_0^t e^{- i (k^3-k_1^3-k_2^3) s}
  v_s(k_1) v_s(k_2) \, ds, \quad t \in [0,T], k \neq 0 .
\end{equation}
 We restrict our attention to
initial conditions such that $v_0(0) = 0$. By calling $\dot X$ the
bilinear operator in the r.h.s.
this equation takes the abstract form
$$
v_t = v_s + \int_s^t \dot X_\sigma(v_\sigma,v_\sigma) d\sigma, \qquad
t,s \in [0,T].
$$
where the paths take values in the Hilbert space $H_\alpha$ of complex Fourier
coefficients $\varphi(k)$ with $\varphi(0)=0$ and $\varphi(-k) =
\overline{\varphi(k)}$ endowed with the scalar product
$
\langle \varphi_1, \varphi_2 \rangle_\alpha = \sum_{k \neq 0}
k^{2\alpha} \varphi_1(-k) \varphi_2(k)
$.

By iteratively substituting the unknown in this integral equation we
obtain an expansion whose first terms looks like
\begin{equation}
\label{eq:tree-series-kdv}
  \begin{split}
v_t & = v_s + \int_s^t d\sigma \dot X_\sigma(v_s,v_s)
+ 2 \int_s^t d\sigma \dot X_\sigma(v_s,\int_s^\sigma d\sigma_1 \dot
X_{\sigma_1}(v_s,v_s))
\\ &\qquad +\int_s^t d\sigma \dot X_\sigma(\int_s^\sigma d\sigma_1 \dot X_{\sigma_1}(v_s,v_s),\int_s^\sigma d\sigma_2 \dot X_{\sigma_2}(v_s,v_s))    
\\ &\qquad +4 \int_s^t d\sigma \dot X_\sigma(v_s,\int_s^\sigma
d\sigma_1 \dot X_{\sigma_1}(v_s,\int_s^{\sigma_1} d\sigma_2 \dot
X_{\sigma_2}(v_s,v_s))  + r_{ts}  
  \end{split}
\end{equation}
where $r_{ts}$ stands for the remaining terms in
the expansion. Denote with $\BT$ the set of (unlabeled) planar rooted
trees with at most two branches at each node. A planar tree is a
rooted tree endowed with an ordering of the branches at each node. 
Then each of the terms in this expansion can be associated
to a tree in $\BT$  and we can define recursively
multi-linear operators $X^{\tau}$ as
$$
X^{\troot}_{ts}(\varphi_1,\varphi_2) = \int_s^t \dot
X_\sigma(\varphi_1,\varphi_2) d\sigma;
$$ 
$$
X^{[\tau^1]}_{ts}(\varphi_1,\dots,\varphi_{m+1}) = \int_s^t \dot
X_\sigma(X^{\tau^1}_{\sigma
  s}(\varphi_1,\dots,\varphi_m),\varphi_{m+1})
 d\sigma
$$
and
$$
X^{[\tau^1 \tau^2]}_{ts}(\varphi_1,\dots,\varphi_{m+n}) = \int_s^t \dot
X_\sigma(X^{\tau^1}_{\sigma s}(\varphi_1,\dots,\varphi_m),X^{\tau^2}_{\sigma s}(\varphi_{m+1},\dots,\varphi_{m+n})) d\sigma.
$$
Eq.~(\ref{eq:tree-series-kdv}) has the form
\smalltrees
\begin{equation}
  \label{eq:kdv-incr-1}
\der v_{ts} = X^{\tsroot}(v^{\times 2})_{ts} + X^{\aabb}(v^{\times 3})_{ts} +
X^{\aaabbb}(v^{\times 4})_{ts} +  X^{\aababb}(v^{\times 4})_{ts} +r_{ts}  
\end{equation}
as an equation in $\CC_2$ where increments take values in $H_\alpha$ and where we let $v^{\times n}_s =
(v_s,\dots,v_s)$ ($n$ times). The operators $X^\tau$ satisfy multiplicative relations 
$$
\der X^{\aabb}(\varphi_1,\varphi_2,\varphi_3) = X^{\tsroot}(X^{\tsroot}(\varphi_1,\varphi_2),\varphi_3) ,
$$
$$
\der X^{\aaabbb}(\varphi_1,\varphi_2,\varphi_3,\varphi_4) =
X^{\tsroot}(X^{\aabb}(\varphi_1,\varphi_2,\varphi_3),\varphi_4) + X^{\aabb}(X^{\tsroot}(\varphi_1,\varphi_2),\varphi_3,\varphi_4) ,
$$
and
\begin{equation*}
  \begin{split}
\der X^{\aababb}& (\varphi_1,\varphi_2,\varphi_3,\varphi_4) =
X^{\tsroot}(X^{\tsroot}(\varphi_1,\varphi_2),X^{\tsroot}(\varphi_3,\varphi_4))
\\ & \qquad +
X^{\aabb}(\varphi_1,\varphi_2,X^{\tsroot}(\varphi_3,\varphi_4)) + 
X^{\aabb}(\varphi_3,\varphi_4,X^{\tsroot}(\varphi_1,\varphi_2))    
  \end{split}
\end{equation*}
where we used the symmetry of the operator $\dot X$ to obtain this
last equation. These relations have much in common with the analogous
relations for branched rough paths, however here the additional
information of the position of the various arguments must be taken
into account in the combinatorics of the reduced coproduct. It would
be interesting to determine  a Hopf algebra structure  on $\BT$
which could account this additional information in a general way. 

Note that $\dot X$ is unbounded on $H_\alpha$ while
it is possible to prove that $X^{\troot}$ and $X^{\aabb}$ are bounded,
in particular for the simplest of them, $X^\troot$, we have the bound
$$
|X^{\troot}_{ts}(\varphi_1,\varphi_2)|_\alpha \le C |t-s|^\gamma
|\varphi_1|_\alpha |\varphi_2|_\alpha 
$$
with the parameters $\gamma$ and $\alpha$ satisfying
\begin{equation}
  \label{eq:range}
\text{$ \gamma < 1/2$, $\alpha \in \RR$ with $\gamma-\alpha< 1$ and
$3 \gamma-\alpha<3/2$.}
\end{equation}
Moreover we are able to decompose the second member of the hierarchy
$X^{\aabb}$ as 
\begin{equation}
\label{eq:decomp-kdv-2}
X^{\aabb}(\varphi_1,\varphi_2,\varphi_3)  = {\widehat X}^{\aabb}(\varphi_1,\varphi_2,\varphi_3)+\Phi(\varphi_1,\varphi_2)\varphi_3+\Phi(\varphi_1,\varphi_3)\varphi_2
\end{equation}
where for the same range of $\gamma,\alpha$ as in~(\ref{eq:range}) we have
$$
|\hat X^{\aabb}_{ts}(\varphi_1,\varphi_2,\varphi_3)|_\alpha \le C |t-s|^{2\gamma}
|\varphi_1|_\alpha |\varphi_2|_\alpha |\varphi_3|_\alpha 
$$
while the operator $\Phi_{ts} : H_\alpha \times H_\alpha \to \CC$ is bounded only for $\alpha \ge -1/2$ and
$$
|\Phi_{ts}(\varphi_1,\varphi_2)| \le C |t-s| |\varphi_1|_\alpha |\varphi_2|_\alpha .
$$
The analysis of the higher order operators has not yet been
performed. However already at this stage something can be said if we
take $3\gamma > 1$ since we
 are naturally led to consider
eq.~(\ref{eq:kdv-incr-1}) as a increment equation and rewrite it
using the sewing map and  the operators $X$ (only up to second order) obtaining the equation
\begin{equation}
  \label{eq:lambda-kdv}
\der v = (1-\Lambda \der)[X^{\tsroot}(v^{\times 2}) + X^{\aabb}(v^{\times 3})]  
\end{equation}
which  can be solved by fixed point methods in $H_\alpha$ for any
$\alpha > -1/2$ (cfr.~(\ref{eq:range}) and the condition on $\Phi$).

The condition $\alpha \ge -1/2$ is essentially imposed by the operator
$\Phi$ appearing in the decomposition~(\ref{eq:decomp-kdv-2}) of
$X^{\aabb}$. This constraint it is linked with a resonance phenomenon
which apperars in the scattering by the non-linear term involving
four waves and which hint to the  fact that the KdV equation
is not uniformly wellposed in $H_\alpha$ for $\alpha < -1/2$~\cite{MR2018661,MR1969209}.

We have replaced the differential and integral approach to the study of this equation by an approach based on an operator valued rough path and the increment complex. It is interesting then to look how the properties of the dynamical system reflect in this unusual approach. As an example let us consider conservation laws.

The KdV equation formally conserves the $H_0$ norm. This conservation
law imposes additional algebraic relations to the operators $X$:
it is not difficult to prove that we have
\begin{equation}
  \label{eq:conservation1}
\langle \varphi_1,  X^\troot_{ts}(\varphi_2,\varphi_3)\rangle_0 
+\langle \varphi_2, X^\troot_{ts}(\varphi_1,\varphi_3)\rangle_0 
+\langle \varphi_3,  X^\troot_{ts}(\varphi_2,\varphi_1)\rangle_0 
= 0  
\end{equation}
and that
\begin{equation}
  \label{eq:conservation2}
2 \langle \varphi, X^2_{ts}(\varphi,\varphi,\varphi)\rangle_0
+   \langle X_{ts}(\varphi,\varphi), X_{ts}(\varphi,\varphi)\rangle_0 = 0  , 
\end{equation}
where all the test functions belong to $H_0$. To see that these two
relations imply the $H_0$ conservation law for solutions we will prove
that $\der \langle v,v \rangle_0  = 0$ when $v$ satisfy
\begin{equation}
  \label{eq:lambda-kdv-r}
\der v = X^{\tsroot}(v^{\times 2}) + X^{\aabb}(v^{\times 3}) + \crr  
\end{equation}
cfr. eq.~(\ref{eq:lambda-kdv}). Let us compute explicitly 
$[\der \langle v,v \rangle_0]_{ts}   = \langle v_t,v_t \rangle_0 - \langle v_s,v_s \rangle_0
=  2 \langle \der v_{ts},v_s \rangle_0 + \langle \der v_{ts},\der v_{ts} \rangle_0    $
.
Substituting in this expression the equation~(\ref{eq:lambda-kdv-r}) we get
\begin{equation*}
  \begin{split}
[\der \langle v,v \rangle_0]_{ts}
& = 2  \langle X^\troot_{ts}(v_s,v_s)+X^{\aabb}_{ts}(v_s,v_s,v_s)   ,v_s \rangle_0 
\\ & \qquad + \langle X^\troot_{ts}(v_s,v_s),X^\troot_{ts}(v_s,v_s) \rangle_0 
 + \crr  .
  \end{split}
\end{equation*}
The relation~(\ref{eq:conservation1}) 
implies that $\langle v_s,
X^\troot_{ts}(v_s,v_s) \rangle_0 = 0$ while eq.~(\ref{eq:conservation2})
allows to cancel the $X^{\aabb}$ term with the quadratic $X^\troot$ term. After the
cancellations the increment of the $H_0$ norm squared is  then $[\der
\langle v,v \rangle_0] \in \CC_2^{1+} $ but this means that it must be
zero and that $\langle v_t,v_t \rangle_0 = \langle v_0,v_0 \rangle_0$
for any $t \ge 0$.

\subsection{Navier-Stokes-like equations}

We consider the NS equation in $\RR^3$ which going in Fourier space
can be written
\begin{equation}
  \label{eq:NS-basic}
v_t(k) = e^{-|k|^2 t} v_0(k) + i \int_0^t e^{-|k|^2 (t-s)} \int_{\RR^3}dk' \langle k, v_s(k-k')\rangle P_k v_s(k') \,ds
\end{equation}
where $v_t$ is the Fourier transform of
the velocity field, $\langle \cdot ,\cdot \rangle$ is the scalar
product in $\mathbb{C}^3$ and $P_k: \CC^3 \to \CC^3$ is the projection
on the directions orthogonal to the vector $k \in \RR^3$, i.e.
$
P_k a = a - \langle k,a\rangle k |k|^{-1}
$.
Eq.~(\ref{eq:NS-basic}) will be studied in the spaces
$\Phi(\alpha)$, $\alpha \in [2,3)$ where $v \in \Phi(\alpha)$ if $v \in C(\RR^3;\CC^3)$ with $k \cdot v(k) = 0$ and 
$\|v\|_\alpha = \sup_{k\in\RR^3} |k|^\alpha |v(k)| < \infty $.
We will write $\alpha = 2+\eps$ with $\eps \in [0,1)$.
The spaces $\Phi(\alpha)$ can contain solutions with infinite energy
and enstrophy so classical results about existence and uniqueness
do not apply.
 Sinai~\cite{Sinai1,Sinai2,Sinai3}, studied
eq.~(\ref{eq:NS-basic}) in $\Phi(\alpha)$ with $\alpha > 2$, showing
that there is existence of unique local
solutions and that these solutions survive for arbitrary large time if
the initial condition is small enough.
Related works on NS are those of Le Jan and Sznitman~\cite{LS},
Cannone and Planchon~\cite{CP} and  the reviews of Bhattacharya et
al. in~\cite{MR1997593} and Waymire~\cite{MR2121794}.

Following~\cite{MR2227041} we will describe the  representation for these solutions as
series indexed by planar binary trees. The use of trees to rigoroulsy analyze the
NS equation has been somewhat pioneered by Gallavotti~\cite{Gall}.

\bigskip
The NS equation can be cast in the abstract form
\begin{equation}
  \label{eq:NS-c-abstract}
u_t = S_t u_0 +  \int_0^t S_{t-s}  B(u_s,u_s)\,ds.
\end{equation}
where $S$ is a bounded semi-group on $\Phi(\alpha)$ and $B$ is a
symmetric bilinear operator which is usually defined only on a
subspace of $\Phi(\alpha)$. 
 Here we cannot proceed as in the KdV case by going to the interaction
 picture since $S$ is only a semi-group, so we must cope with
the convolution directly. 
In~\cite{MR2227041} we showed that the solutions of this equation in the case of the 3d NS equation have
the norm convergent series representation
\begin{equation}
  \label{eq:ns-series}
u_t = S_t u_0 + \sum_{\tau \in \BT} X^{\tau}_{t0}(u_0^{\times \theta(\tau)})   
\end{equation}
where $\theta(\tau)$ is a degree function defined by $
\theta(\troot) = 2$, $\theta([\tau]) = 1+ \theta(\tau)$, $
\theta([\tau_1 \tau_2]) = \theta(\tau^1) + \theta(\tau^2)  
$ and the $\theta(\tau)$-multilinear operators
$X^{\tau}$ have recursive definition
$$
X^{\tsroot}_{ts} (\varphi^{\times 2}) = \int_s^t S_{t-u}
B(S_{u-s}\varphi,S_{u-s}\varphi) du
$$
$$
X^{[\tau^1]}_{ts}(\varphi^{\times (\theta(\tau^1)+1)}) = \int_s^t
S_{t-u}B(X^{\tau^1}_{us}(\varphi^{\times \theta(\tau^1)}),\varphi) du
$$
and
$$
X^{[\tau^1 \tau^2]}_{ts}(\varphi^{\times (\theta(\tau^1)+\theta(\tau^2))}) = \int_s^t
S_{t-u}B(X^{\tau^1}_{us}(\varphi^{\times \theta(\tau^1)}),X^{\tau^2}_{us}(\varphi^{\times \theta(\tau^2)})) du
$$
and by induction we can prove that, for any $\eps \in [0,1)$ these
operators are bounded by
\begin{equation}
\label{eq:bounds-NS}
|X^\tau_{ts}(h^{\theta(\tau)})(k)| \le C_{\tau} \frac{e^{-|k|^2 (t-s) / (|\tau|+1)}}{|k|^\alpha}
(t-s)^{|\tau| \eps /2} \|h\|_\alpha^{\theta(\tau)}  
\end{equation}
where the constants $C_{\tau}$ can be chosen as 
\begin{equation}
  \label{eq:asymp-q-ns}
C_\tau = A^{|\tau|} (\tau!)^{-\eps/2}.  
\end{equation}
for some other constant $A > 0$ depending only on $\eps$.

Due to the presence of the convolution integral these $X$ operators
does not behaves nicely with respect to
the coboundary $\der$. In~\cite{TindelGubinelli} we introduced cochain complex $(\hat
C_*,\tilde \der)$ adapted to the study of such convolution integrals
where the coboundary $\tilde \der$ is obtained from $\der$ by a "twisting" involving the semigroup.
 There exists also a corresponding convolutional sewing map $\tilde
\Lambda$ which
provide an appropriate inverse to $\tilde \der$. Algebraic relations for
these iterated integrals have then  by-now familiar expressions, e.g.:
$$
\tilde \der X^{\aabb}(\varphi^{\times 3}) =
X^{\tsroot}(X^{\tsroot}(\varphi^{\times 2}), S\varphi) 
$$
and so on. Indeed $X$ can be considered as a branched rough path and
the asymptotic behaviour~(\ref{eq:asymp-q-ns}) supports somewhat the
conjectured asymptotics~(\ref{eq:q-conj}).

\bigskip

The series representation~(\ref{eq:ns-series}) together with the
bounds~(\ref{eq:bounds-NS}) imply that
\begin{equation}
  \label{eq:bounds-sol-NS}
|u_t(k)| 
\le  |S_t u_0 (k)| + \sum_{\tau \in \BT}\frac{
  B^{|\tau|}}{\sigma(\tau) \tau!} \frac{e^{-(k^2 t)/(|\tau|+1)}}{|k|^\alpha} t^{\eps
  |\tau|/2} \|u_0\|_\alpha^{\theta(\tau)}   
\end{equation}
By induction we can prove that $\gamma(\tau) \ge 2^{|\tau|-1}$ where
equality holds for the binary trees for which every path from the root
to the leaves has the same length. This estimate together with
\begin{equation}
  \label{eq:theta-bound}
 (|\tau|+1)/2 \le
\theta(\tau) \le
|\tau|+1  
\end{equation}
easily proven by induction on $|\tau|$ give
\begin{equation}
  \label{eq:bounds-sol-NS-2}
|u_t(k)| 
\le  |S_t u_0 (k)| + \sum_{n \ge 1}
Z_n  B^{n}\frac{ e^{-(k^2 t)/(n+1)}}{|k|^\alpha} t^{\eps
 n/2} \|u_0\|_\alpha^{(n+1)/2} (1+\|u_0\|_\alpha)^{(n+1)/2}   
\end{equation}
for some different constant $B$ and where $Z_n$ is the number of trees
in $\BT$ with $n$ vertices for which we have the estimate
$
Z_n \le D^n (n+1)^{-3/2}
$ for a constant $D>0$. So the series~(\ref{eq:ns-series}) is
controlled by the geometric series~(\ref{eq:bounds-sol-NS-2}) and converges
in norm $t$ is small or $\|h\|_\alpha$ is small. In the case $\eps =
0$ the dependence in time of the r.h.s. is bounded and so the series
converges for all time if the initial condition is small enough.
The series gives also additional informations on the global solution
when $\eps =0$:
\begin{itemize}
\item[a)]  \emph{dissipation}: for fixed $k \in \RR^3\backslash \{0\}$, $
\lim_{t \to \infty } |u_t(k)| = 0  $;
\item[b)] \emph{smoothness}: for fixed $t > 0$, there exists two constants $C_3,C_4$ such that
 $|v_t(k)| \le C_3 e^{-C_4 |k|\sqrt{t}}$ as $|k| \to \infty$.
\end{itemize}
In particular this second property can be proved by the Laplace method applied to the majorizing series\ref{eq:bounds-sol-NS-2}.

\medskip
Another interesting way to analyze the NS series~(\ref{eq:ns-series}) is to note that different classes of trees give different
contributions. We define \emph{simple} trees the trees with at most
one branch at each vertex, i.e. of the form $[\cdots [\bullet]\cdots
]$. \emph{Short} trees are instead trees for which at  each  vertex we
have two branches, each of which carries (asymptotically) a fixed proportion ($\alpha$ or $1-\alpha$) of the
vertice and without loosing generality we consider  $\alpha \in (0,1/2)$.
We will denote $\BT_0$ the set of simple trees and $\BT_\alpha$ the
set  of short trees corresponding to the proportion $\alpha$. 
The distinction between these classes of trees is relevant when
discussing the asymptotic behavior of the tree factorial. Indeed
for $\tau \in \BT_0$ we have $\tau! = |\tau|!$ while for any $\alpha
\in (0,1/2)$ there exists constants $D_1,D_2,D_3,D_4$ such that, for any
$\tau \in \BT_\alpha$ we have
$$
 D_3   |\tau|^{-1} D_4^{|\tau|} \le \gamma(\tau) \le D_1   |\tau|^{-1} D_2^{|\tau|}.
$$  
This different behavior  is responsible for
different convergence properties of the sum~(\ref{eq:ns-series}) when
restricted to simple or short trees. Indeed it is possible to prove
that the series restricted to simple trees is convergent for all times
whatever the size of the initial condition while the estimate for
short trees do not ensure this important property. In some sense the
difficulty of finding global solution of NS is due to the presence of
arbitrarily large short trees in the expansion.
A similar phenomenon is observed in~\cite{Sinai3} and exploited in~\cite{MR2390325} to prove blow-up for the complex Navier-Stokes equation by a renormalization group argument.


\begin{thebibliography}{10}

\bibitem{MR0454968}
Kuo~Tsai Chen.
\newblock Iterated path integrals.
\newblock {\em Bull. Amer. Math. Soc.}, 83(5):831--879, 1977.

\bibitem{MR1654527}
Terry~J. Lyons.
\newblock Differential equations driven by rough signals.
\newblock {\em Rev. Mat. Iberoamericana}, 14(2):215--310, 1998.

\bibitem{MR2036784}
Terry Lyons and Zhongmin Qian.
\newblock {\em System control and rough paths}.
\newblock Oxford Mathematical Monographs. Oxford University Press, Oxford,
  2002.
\newblock Oxford Science Publications.

\bibitem{MR2053040}
Antoine Lejay.
\newblock An introduction to rough paths.
\newblock In {\em S\'eminaire de Probabilit\'es XXXVII}, volume 1832 of {\em
  Lecture Notes in Math.}, pages 1--59. Springer, Berlin, 2003.

\bibitem{ramif}
M.~Gubinelli.
\newblock Ramification of rough paths.
\newblock {\em J. Diff. Eq.}, 2008.
\newblock to appear.

\bibitem{MR2227041}
Massimiliano Gubinelli.
\newblock Rooted trees for 3{D} {N}avier-{S}tokes equation.
\newblock {\em Dyn. Partial Differ. Equ.}, 3(2):161--172, 2006.

\bibitem{kdv}
M.~Gubinelli.
\newblock Rough solutions of the periodic {K}orteweg-de {V}ries equation.
\newblock 2006.

\bibitem{MR2091358}
M.~Gubinelli.
\newblock Controlling rough paths.
\newblock {\em J. Funct. Anal.}, 216(1):86--140, 2004.

\bibitem{MR0305608}
J.~C. Butcher.
\newblock An algebraic theory of integration methods.
\newblock {\em Math. Comp.}, 26:79--106, 1972.

\bibitem{MR2106008}
Ch. Brouder.
\newblock Trees, renormalization and differential equations.
\newblock {\em BIT}, 44(3):425--438, 2004.

\bibitem{MR1660199}
Alain Connes and Dirk Kreimer.
\newblock Hopf algebras, renormalization and noncommutative geometry.
\newblock {\em Comm. Math. Phys.}, 199(1):203--242, 1998.

\bibitem{TindelGubinelli}
M.~Gubinelli and S.~Tindel.
\newblock Rough evolution equations.
\newblock 2006.

\bibitem{MR2018661}
Michael Christ, James Colliander, and Terrence Tao.
\newblock Asymptotics, frequency modulation, and low regularity ill-posedness
  for canonical defocusing equations.
\newblock {\em Amer. J. Math.}, 125(6):1235--1293, 2003.

\bibitem{MR1969209}
J.~Colliander, M.~Keel, G.~Staffilani, H.~Takaoka, and T.~Tao.
\newblock Sharp global well-posedness for {K}d{V} and modified {K}d{V} on
  {$\Bbb R$} and {$\Bbb T$}.
\newblock {\em J. Amer. Math. Soc.}, 16(3):705--749 (electronic), 2003.

\bibitem{Sinai1}
Ya.~G. Sinai.
\newblock On local and global existence and uniqueness of solutions of the 3{D}
  {N}avier-{S}tokes system on {$\Bbb R\sp 3$}.
\newblock In {\em Perspectives in analysis}, volume~27 of {\em Math. Phys.
  Stud.}, pages 269--281. Springer, Berlin, 2005.

\bibitem{Sinai2}
Yakov Sinai.
\newblock Power series for solutions of the {$3D$}-{N}avier-{S}tokes system on
  {${\bf R}\sp 3$}.
\newblock {\em J. Stat. Phys.}, 121(5-6):779--803, 2005.

\bibitem{Sinai3}
Ya.~G. Sina{\u\i}.
\newblock A diagrammatic approach to the 3{D} {N}avier-{S}tokes system.
\newblock {\em Uspekhi Mat. Nauk}, 60(5(365)):47--70, 2005.

\bibitem{LS}
Y.~Le~Jan and A.~S. Sznitman.
\newblock Stochastic cascades and {$3$}-dimensional {N}avier-{S}tokes
  equations.
\newblock {\em Probab. Theory Related Fields}, 109(3):343--366, 1997.

\bibitem{CP}
Marco Cannone and Fabrice Planchon.
\newblock On the regularity of the bilinear term for solutions to the
  incompressible {N}avier-{S}tokes equations.
\newblock {\em Rev. Mat. Iberoamericana}, 16(1):1--16, 2000.

\bibitem{MR1997593}
Rabi~N. Bhattacharya, Larry Chen, Scott Dobson, Ronald~B. Guenther, Chris Orum,
  Mina Ossiander, Enrique Thomann, and Edward~C. Waymire.
\newblock Majorizing kernels and stochastic cascades with applications to
  incompressible {N}avier-{S}tokes equations.
\newblock {\em Trans. Amer. Math. Soc.}, 355(12):5003--5040 (electronic), 2003.

\bibitem{MR2121794}
Edward~C. Waymire.
\newblock Probability \& incompressible {N}avier-{S}tokes equations: an
  overview of some recent developments.
\newblock {\em Probab. Surv.}, 2:1--32 (electronic), 2005.

\bibitem{Gall}
Giovanni Gallavotti.
\newblock {\em Foundations of fluid dynamics}.
\newblock Texts and Monographs in Physics. Springer-Verlag, Berlin, 2002.
\newblock Translated from the Italian.

\bibitem{MR2390325}
Dong Li and Ya.~G. Sinai.
\newblock Blow ups of complex solutions of the 3{D} {N}avier-{S}tokes system
  and renormalization group method.
\newblock {\em J. Eur. Math. Soc. (JEMS)}, 10(2):267--313, 2008.

\end{thebibliography}

\end{document}